\documentstyle{article}
\def\bz{{\bf Z}}
\def\tp{{\widetilde{\phi}}}
\begin{document}
\vskip1cm
\begin{center}
{\Large{Determinant formulas for the $\tau$-functions}
\vskip2mm
{of the Painlev\'e equations of type $A$}}
\vskip1cm
{Yasuhiko Yamada}
\vskip5mm
{Department of Mathematics, Kobe University, Rokko, Kobe 657-8501, Japan}
\end{center}
\vskip1cm
{\it{Abstract}}.
Explicit determinant formulas are presented for the
$\tau$-functions of  the generalized Painlev\'e equations of type $A$.
This result allows an interpretation of the $\tau$-functions
as the Pl\"ucker coordinates of the universal Grassmann manifold.
\vskip1cm

\section{Introduction}

For each generalized Cartan matrix of affine type $A=(a_{ij})_{ij \in I}$,
we introduced a representation of the Weyl group $W(A)$ 
on the rational functions of variables $\alpha_i, f_i, \tau_i$, $(i \in I)$
\cite{affine}. 
The representation is characterized by the action of the generator
$s_i$, $(i \in I)$, such that
$$
s_i(\alpha_j)=\alpha_j-\alpha_i a_{ij}, \quad
s_i(f_j)=f_j+{\alpha_i \over f_i}u_{ij}, \quad
s_i(\tau_j)=\tau_j (f_j \prod_{k \in I} \tau_k^{-a_{kj}})^{\delta_{ij}}.
$$
under certain conditions on the constants $(u_{ij})_{ij \in I}$.

This representation is a generalization of the
B\"acklund transformations of the Painlev\'e equations 
$P_{IV},P_{V}$ and $P_{VI}$.
As to the root systems of type $A^{(1)}_l$,
it is also known that there exist Painlev\'e type
differential (or difference) system which has the $W(A)$ symmetry 
\cite{affine}\cite{higher}.
In the context of the Painlev\'e equations, 
the variables $\alpha_i$, $f_i$ and $\tau_i$ play
the role of parameters (or the discrete time variables), 
dependent variables and the $\tau$-functions, respectively.  
In \cite{affine}, it is conjectured that the $\tau$-functions 
have strong regularity.
This regularity is crucial for the differential (or difference) 
systems, since it should  be closely related to the Painlev\'e 
(or the singularity confinement) properties.

In this paper, we prove the regularity conjecture in the case of
affine Weyl groups of type $A^{(1)}_l$ and $A_{\infty}$,
by constructing explicit determinant formulas for the $\tau$-functions.
The determinant formulas allow an interpretation of the $\tau$-functions
as Pl\"ucker coordinates of the universal Grassmann manifold 
as in the theory of KP hierarchy \cite{KP}.
\section{The $\tau$-functions}

We formulate our results in the case of $A_{\infty}$.
For the reduction to $A^{(1)}_l$ case, see Remark 1 in the last section.

Let $K={\bf{C}}(\alpha;f;\tau)$ be the field of rational functions in  
infinitely many variables 
$\alpha=(\alpha_i)$, $f=(f_i)$ and $\tau=(\tau_i)$ where $i \in \bz$.
For each $i \in \bz$, define an automorphism 
$s_i$ of the field $K$ by
$$
s_i(\alpha_i)=-\alpha_i, \quad
s_i(\alpha_{i \pm 1})=\alpha_{i \pm 1} + \alpha_i, \quad
s_i(\alpha_j)=\alpha_j, \ (j \neq i, i \pm 1),
$$
$$
s_i(f_i)=f_i, \quad
s_i(f_{i \pm 1})=f_{i \pm 1} \pm {\alpha_i \over f_i}, \quad
s_i(f_j)=f_j, \ (j \neq i, i \pm 1),
$$
and
$$
s_i(\tau_i)=f_i {\tau_{i-1} \tau_{i+1} \over \tau_i}, \quad
s_i(\tau_j)=\tau_j, \ (j \neq i).
$$
These automorphisms satisfy the relations
$$
s_i^2=1, \quad (s_i s_{i+1})^3=1, \quad s_is_j=s_js_i, \ (j \neq i\pm 1),
$$
and define a representation of the affine Weyl group $W=W(A_{\infty})$ of 
type $A_{\infty}$.
One can extend this representation $W$ to the extended affine Weyl group 
$\widetilde{W}$
by adding the automorphism $\pi$ defined as
$$
\pi(\alpha_i)=\alpha_{i+1}, \quad
\pi(f_i)=f_{i+1}, \quad
\pi(\tau_i)=\tau_{i+1}.
$$
The diagram shift $\pi$ satisfies the relation $\pi s_i=s_{i+1} \pi$.

Let $\Lambda_i$ $(i \in \bz)$ be
the fundamental weights of
$A_{\infty}$ on which the Weyl group ${\widetilde W}$ acts as
$$
s_i(\Lambda_i)=\Lambda_{i+1}+\Lambda_{i-1}-\Lambda_{i}, \quad
s_i(\Lambda_j)=\Lambda_j, \ (j \neq i), \quad
\pi(\Lambda_i)=\Lambda_{i+1}.
$$
We also use the notation $v_i=\Lambda_i-\Lambda_{i-1}$.
Then
$\alpha_i=-\Lambda_{i-1}+2 \Lambda_i-\Lambda_{i-1}=v_i-v_{i+1}$,
and we can put formally
$\Lambda_i=\sum_{j \leq i} v_j$.
The action of $s_i$ on $v_j$ is given by the permutation 
$$
s_i(v_i)=v_{i+1}, \quad
s_i(v_{i+1})=v_{i}, \quad
s_i(v_j)=v_j, \ (j \neq i,i+1).
$$
Furthermore, the Weyl group $W$ is identified with the
infinite symmetric group $S_{\infty}$ which permutes $v_j$.

There exist a family of rational functions
$\phi_{w}(\Lambda_j) \in {\bf C}(\alpha; f)$ for 
$w \in {\widetilde W}$ and $j \in \bz$ such that
$$
w(\tau_j)=\phi_{w}(\Lambda_j) \prod_{i \in \bz}\tau_i^{m_i},
$$
where $m_i=\langle \alpha_i,w \Lambda_j \rangle \in \{0,\pm1 \}$ 
is the coefficient of $\Lambda_i$ in $w \Lambda_j$.
The conjecture in \cite{affine} states that the function $\phi_w(\Lambda_j)$ is
a polynomial in $\alpha$ and $f$ with integral coefficients. 
We will prove this conjecture by constructing an explicit
determinant formula for $\phi_w(\Lambda_j)$. 
Since any function $\phi_w(\Lambda_j)$ is easily obtained from the
case $j=0$ by the shift $\phi_w(\Lambda_j)=\pi^j\big(\phi_{\pi^{-j} w
\pi^{j}}(\Lambda_0)\big)$, we will concentrate on the functions $\phi_w(\Lambda_0)$.

\vskip5mm
\noindent
{\bf{Lemma 1.}}
If $w_1 \Lambda_0=w_2 \Lambda_0$ for $w_1,w_2 \in W$, then
$w_1(\tau_0)=w_2(\tau_0)$.
\vskip5mm
\noindent
{\it{Proof.}}
The condition $w_1(\Lambda_0)=w_2(\Lambda_0)$ means that the element
$w=w_2^{-1} w_1 \in W$ is in the stabilizer $W_0$ of $\Lambda_0$.
An element in $W$ belongs to $W_0$ if and only if 
the corresponding permutation of $\{ v_i \ \vert \  i \in \bz \} \simeq \bz$ 
preserves the subset 
$\bz_{\leq0}$.  
Then $W_0$ is a product of permutations of
$\bz_{\leq 0}$ and $\bz_{\geq 1}$.
Hence $w \in W_0=\langle s_j, \ (j \in \bz_{< 0}) \rangle \times
\langle s_j , \ (j \in \bz_{\geq 1}) \rangle$ and $w(\tau_0)=\tau_0$.
$\Box$
\vskip5mm

By this Lemma, the functions $\phi_w(\Lambda_0)$ are parameterized by 
$w \Lambda_0$. We will prepare some notations to describe the 
orbit $W.\Lambda_0=\{w \Lambda_0 \ \vert \ w \in W \}$
whose elements will be parameterized by Young diagrams \cite{ma}.

Let $\lambda=(\lambda_1 \geq \lambda_2 \geq \ldots \geq \lambda_l >0)$ 
be a partition of length $l=l(\lambda)$.
The corresponding Young diagram $Y$ is defined by
$$
Y=\{s=(i,j) \ \vert \ 1 \leq i \leq l, \ 1 \leq j \leq \lambda_i \}.
$$
The transposition $\lambda'=Y'$ is defined by $(i,j) \in Y'$ 
if and only if $(j,i) \in Y$.
For a partition $\lambda$,
the corresponding Frobenius symbol 
$$
(I,J)=(\{i_1>i_2>\ldots>i_k>0\},\{j_1>j_2>\ldots>j_k\geq 0\})
$$
is defined by
$$
i_n=\lambda_n-n+1, \quad j_n=\lambda'_n-n,
$$ 
where $k=\max\{ n \vert \lambda_n \geq n \}$.
We always identify the three notions, partition $\lambda$, 
Young diagram $Y$ and Frobenius symbol $(I,J)$ by these correspondence.

For $w \in W$, the element $w(\Lambda_0)$ can be parameterized by the partitions
$\lambda=Y=(I,J)$ as follows.

Recall that $\Lambda_0=\sum_{j \leq 0}v_j$.
Then we have
$w(\Lambda_0)=\sum_{i \in M} v_i$, where $M$ is a subset of $\bz$
different from $\bz_{\leq 0}$ only by finite elements.
Such subset $M$ (called Maya diagram) corresponds to a Young diagram
$Y$ with Frobenius symbol $(I,J)$ by the rule 
$$
M \cup (-J)=\bz_{\leq 0} \cup I.
$$

In terms of the Maya diagram $M$, the coefficient $m_i$ of $\Lambda_i$ in 
$w(\Lambda_0)=(\sum_{i \in M}-\sum_{{i+1} \in M})\Lambda_i$ is given by 
\begin{itemize}
\item{$m_i=1$, if $i \in M $ and $i+1 \notin M$},
\item{$m_i=-1$, if $i \notin M$ and $i+1 \in M$},
\item{$m_i=0$, if $i,i+1 \in M$ or $i,i+1 \notin M$}.
\end{itemize}
Translating this into the language of the Young diagrams, we have
\vskip5mm
\noindent
{\bf{Lemma 2.}}
Let $w(\Lambda_0)=\sum_{i \in \bz} m_i \Lambda_i$ and let
$Y$ be the corresponding Young diagram,
then the Young diagram $s_i(Y)$ corresponding to $s_iw(\Lambda_0)$
is given as follows,
\begin{itemize}
\item{$s_i(Y)$ is obtained by adding the node with color $i$, if $m_i=1$,}
\item{$s_i(Y)$ is obtained by removing the node with color $i$, if $m_i=-1$,}
\item{$s_i(Y)=Y$, if $m_i=0$,}
\end{itemize}
where the color $k$ of the $(i,j)$-th node is given by $k=j-i$.
\vskip5mm
\noindent
In summary, we have
\vskip5mm
\noindent
{\bf{Proposition 1.}}
Any element in $W.\Lambda_0$ 
can be obtained from $\Lambda_0$ by the action of 
$s_{i_p}\cdots s_{i_1}$ with 
$m_{i_k}=\langle \alpha_{i_k},s_{i_k-1}\cdots s_{i_1}(\Lambda_0) \rangle=1$,
$(k=1,\ldots,p)$.
Hence, the functions
$\phi_w(\Lambda_0)$ are uniquely determined by the cocycle condition
$$
\phi_{s_i w}(\Lambda_0)=s_i\big(\phi_w(\Lambda_0)\big)f_i, \quad (m_i=1),
$$
with the initial condition $\phi_{\it{id}}(\Lambda_0)=1$.
\vskip5mm

It is convenient to introduce another normalization 
$\tp_w(\Lambda_0)$ defined by
$$
\tp_{w}(\Lambda_0)=
{1 \over N_w} \phi_{w}(\Lambda_0).
$$ 
Here the normalization factor
$N_{w}$ is a polynomial in $\alpha$, 
which is defined, in terms of corresponding
Young diagram $w(\Lambda_0) \leftrightarrow \lambda=Y$, as 
$$
N_{w}=\prod_{s =(i,j) \in Y}h(s,\alpha),
\quad
h((i,j), \alpha)=v_{j-\lambda'_j}-v_{\lambda_i-i+1}.
$$
Note that when specialized to $\alpha_i=1$, $h(s,1)$ is nothing but the
hook-length of the node $s \in Y$.

\vskip5mm
\noindent
{\bf{Lemma 3.}}
The normalization factor $N_w$ satisfies the relation
$$
N_{s_k w}=\alpha_k s_k(N_w), \quad (m_k=1).
$$
{\it{Proof.}}
{}From the Lemma 2 and the condition $m_k=1$ 
the Young diagram $s_k(Y)$ is obtained from
$Y$ by adding one node, say $(i_0,j_0)$-th node, with color $k=j_0-i_0$.
Since $\lambda_{i_0}=j_0-1$ and $\lambda'_{j_0}=i_0-1$,
the hook length $h(s,\alpha)=v_{j-\lambda'_j}-v_{\lambda_i-i+1}$ 
for $s=(i,j) \in Y$ contains $v_k$ or $v_{k+1}$ if and only if
\begin{itemize}
\item{$s=(i_0,*)$, with $h(s,\alpha)=v_{*}-v_{k}$,}
\item{$s=(*,j_0)$, with $h(s,\alpha)=v_{k+1}-v_{*}$.}
\end{itemize}
Under the action of $s_k$, these factors are replaced with the hook length of the
same node $s$ in the new diagram $s_k(Y)$.
Multiplying $s_k(N_w)$ by the extra factor $v_k-v_{k+1}=\alpha_k$ 
corresponding to the added node $(i_0,j_0)$, we get $N_{s_k w}$. $\Box$
\vskip5mm
\noindent
As a Corollary of this Lemma 3 and Proposition 1 we have
\vskip5mm
\noindent
{\bf{Proposition 2.}}
The normalized functions $\tp_w(\Lambda_0)$
are determined by the cocycle condition
$$
\tp_{s_i w}(\Lambda_0)=s_i\left(\tp_w(\Lambda_0)\right){f_i \over \alpha_i},
\quad (m_i=1),
$$
with the initial condition $\tp_{\it{id}}(\Lambda_0)=1$.

\vskip5mm
\noindent
{\bf{Examples 1.}}
For the partitions $(0),(1),(2),(1,1),(2,1)$ and $(2,2)$,
the corresponding normalized functions $\tp_w(\Lambda_0)$ are given
as follows.
$$
a_{12}=\tp_{1}(\Lambda_0)=1 \leftrightarrow (0),
$$
$$
a_{13}=\tp_{s_0}(\Lambda_0)=
{f_0 \over \alpha_0} \leftrightarrow (1),
$$
$$
a_{14}=\tp_{s_1 s_0}(\Lambda_0)=
{f_0 f_1-\alpha_1 \over (\alpha_0+\alpha_1)\alpha_1}
\leftrightarrow (2),
$$
$$
a_{23}=\tp_{s_{-1}s_0}(\Lambda_0)=
{f_{-1}f_0+\alpha_{-1} \over (\alpha_{-1}+\alpha_0)\alpha_{-1}}
\leftrightarrow (1,1),
$$
$$
a_{24}=\tp_{s_{-1}s_1s_0}(\Lambda_0)=
{f_{-1}f_0f_1+\alpha_{-1}f_1-\alpha_1 f_{-1}
\over (\alpha_{-1}+\alpha_0+\alpha_1)\alpha_{-1} \alpha_1}
\leftrightarrow (2,1),
$$
and for the partition $(2,2)$, we have
$$
a_{34}=
\tp_{s_0s_{-1}s_1s_0}(\Lambda_0)=
{f_{-1}f_0^2f_1+\alpha_{-1} f_0 f_1-\alpha_1 f_{-1}f_0+\alpha_0 (\alpha_{-1}+\alpha_0+\alpha_1)
\over (\alpha_{-1}+\alpha_0+\alpha_1)(\alpha_0+\alpha_1)(\alpha_{-1}+\alpha_0)\alpha_0}.
$$

It is interesting to note the relations
\begin{equation}
a_{12} a_{34}-a_{13} a_{24}+a_{14} a_{23}=0,
\end{equation}
\begin{equation}
a_{23}=\pi^{-1}(a_{13})a_{13}-\pi^{-1}(a_{14})a_{12}.
\end{equation}
Each of the equations plays a fundamental role in the proof of 
Theorems 1 or 2 respectively.

\section{The determinant formulas}

For integers $p \geq 1$ and $q \geq 0$, put
$$
X_{p,q}=\det
\left[
\matrix{
f_{-q}&1&&&0 \cr
\beta_{-q+1}&f_{-q+1}&1&& \cr
&\ddots&\ddots&\ddots& \cr
&&\beta_{p-2}&f_{p-2}&1& \cr
0&&&\beta_{p-1}&f_{p-1} \cr
}
\right],
$$
where $\beta_j$ $(-q+1 \leq j \leq p-1)$ is given by
$$
\beta_j=\sum_{i=1}^{p-1}\alpha_i=v_j-v_p, \ (j>0), \quad
\beta_j=-\sum_{i=-q}^{j-1}\alpha_i=v_j-v_{-q+1}, \ (j \leq 0).
$$

\vskip5mm
\noindent
{\bf{Lemma 4.}}
Put $w_{p,q}=(s_{q}\cdots s_{-1})(s_{p-1}\cdots s_0)$, which
corresponds to the hook diagram 
$w_{p,q}(\Lambda_0) \leftrightarrow \lambda=(\{p\},\{q\})$, then
$$
\phi_{p,q}:=\phi_{w_{p,q}}(\Lambda_0)=X_{p,q}.
$$
\vskip5mm
\noindent
{\it{Proof.}}
For $p=1$ and $q=0$, the formula is trivially satisfied, $X_{1,0}=f_0$.
We need to check the following relations for $p \geq 1$ and $q \geq 0$.
\begin{equation}
X_{p+1,q}=s_{p}(X_{p,q})f_p, \quad
X_{p,q+1}=s_{-(q+1)}(X_{p,q}) f_{-(q+1)}.
\end{equation}
Let us prove the first one.
By definition
$$
X_{p+1,q}=\det
\left[
\matrix{
\ddots&\ddots&\cr
\ddots&f_{p-2}&1&\cr
&\alpha_{p-1}+\alpha_p&f_{p-1}&1&\cr
&&\alpha_p&f_{p}
}
\right],
$$
expanding this with respect to the last row, we have
$$
=\det
\left[
\matrix{
\ddots&\ddots&\cr
\ddots&f_{p-2}&1&\cr
&s_p(\alpha_{p-1})&f_{p-1}
}
\right] f_p-
\det
\left[
\matrix{
\ddots&\ddots&\cr
\ddots&f_{p-2}&1&\cr
&s_p(\alpha_{p-1})&1
}
\right] \alpha_p,
$$
which is nothing but
$$
s_p(X_{p,q})=f_p
\det
\left[
\matrix{
\ddots&\ddots&\cr
\ddots&f_{p-2}&1&\cr
&s_p(\alpha_{p-1})&f_{p-1}-{\alpha_p \over f_p}
}
\right].
$$
Thus the first relation is proved.
The second one is similar.$\Box$
\vskip5mm
\noindent
Using the notations above, we can state our main result as follows.
\vskip5mm
\noindent
{\bf{Theorem 1.}}
For any $w \in  W$, the normalized function 
$\tp_w(\Lambda_0)$ is given by the
following determinant
$$
\tp_w(\Lambda_0)=\det \left( \tp_{p,q} \right) _{p \in I,q \in J},
$$
where $(I,J)$ is the Frobenius symbol of $\lambda=Y$ corresponding to 
$w(\Lambda_0)$.
\vskip5mm
\noindent
Equivalently, we also have the following Jacobi-Trudi type (\cite{ma}) 
formulas
\vskip5mm
\noindent
{\bf{Theorem 2.}}
For any $w \in  W$, the normalized function $\tp_w(\Lambda_0)$ 
is given by the following determinant
$$
\tp_w(\Lambda_0)=\det 
\left( h_{\lambda_i-i+j}^{(1-j)} \right)_{1 \leq i,j \leq l(\lambda)},
$$
where $\lambda$ is the partition corresponds to $w(\Lambda_0)$ and
$h_k^{(j)}=\pi^j(\tp_{k,0})$ is the normalized function for the single row
$\lambda=(k)$.
\vskip5mm
\noindent
The proofs of these theorems are given in the next section.
\vskip5mm
\noindent
{\bf{Corollary.}}
For any $w \in W$, the function $\phi_w(\Lambda_0)$ is
a polynomial in $\alpha$ and $f$ with integral coefficients.
The leading term with respect to $f$ is $\prod_{i \in \bz} f_i^{\nu_i}$,
where $\nu_i$ is the number of nodes with color $i$ in 
$Y \leftrightarrow w(\Lambda_0)$.
\vskip5mm
\noindent
{\it{Proof.}}
We will prove by induction on the length of $w \in W$.
$\phi_{\it{id}}(\Lambda_0)=1$ is in $\bz[\alpha ; f]$.
Assume that $R(\alpha;f)=\phi_{w}(\Lambda_0) \in \bz[\alpha ; f]$.
By Proposition 1, the function 
$S(\alpha;f)=\phi_{s_i w}(\Lambda_0)$ for $m_i=1$ is given by
$$
S(\ldots,\alpha_{i-1}+\alpha_{i},-\alpha_i,\alpha_{i+1}+\alpha_i,\ldots;
\ldots,f_{i-1}-{\alpha_i \over f_i},f_i,f_{i+1}+{\alpha_i \over f_i},\ldots) f_i,
$$
which belongs to $\bz[\alpha; f; {1 \over f_i}]$.
The condition on the leading term of $S$ follows from that on $R$. 
On the other hand, from Theorem 1, we see that 
$\phi_{s_i w}(\Lambda_0)=N_{s_i w} \tp_{s_i w}(\Lambda_0)$ is a polynomial in $f$, 
hence the function $S(\alpha;f)$ also belongs to $\bz[\alpha; f]$. 
$\Box$

\section{Proof of the Theorems}
\noindent
{\bf{Lemma 5.}}
We have
$$
z^{n-1}
\det
\left[
\matrix{
(a_{ij})_{1 \leq i,j \leq n}&(x_i)_{1 \leq i \leq n} \cr
(y_j)_{1 \leq j \leq n}&z \cr}
\right]=
\det
\left[ 
\big(a_{ij} z-x_i y_j \big)_{1 \leq i,j \leq n}
\right].
$$
\vskip5mm
\noindent
{\it{Proof.}}
Let $e_i$, ($0 \leq i \leq n$) be a basis in ${\bf{C}}^{n+1}$.
For $1 \leq i \leq n$, put
$$
a_i=\eta_i+x_i e_0, \quad \eta_i=\sum_{j=1}^n a_{ij} e_j,
\quad
a_0=\xi+z e_0, \quad \xi=\sum_{j=1}^n y_j e_j,
$$
$$
{\hbox{and}} \quad b_i=\sum_{j=1}^n (z a_{ij}-x_i y_j) e_j=z \eta_i-x_i \xi.
$$
Then the both hand sides of the identity are the 
coefficients of $e_0 \wedge \cdots \wedge e_n$ in 
$$
(LHS)=z^{n-1} a_0 \wedge a_1 \wedge \cdots \wedge a_n, \quad
(RHS)=e_0 \wedge b_1 \wedge \cdots \wedge b_n.
$$
It is easy to see that these two are the same and equal to
$$
z^n e_0 \wedge \eta_1 \wedge \cdots \wedge \eta_n-z^{n-1} \sum_{i=1}^n 
e_0 \wedge \eta_1 \wedge \cdots \wedge x_i \xi \wedge \cdots \wedge \eta_n.
$$
$\Box$
\vskip5mm
\noindent
{\it{Proof of the Theorem 1.}}
To prove the theorem, it is enough to check that the determinant satisfy the
transformation properties in Proposition 2.

For the actions of $s_i$, $i \neq 0$, the transformation properties
directly follows from eq.(3) for the hook $X_{p,q}$ in Lemma 4, 
since the $s_i$
acts only on the single row ($i>0$) or single column ($i<0$).

For the $s_0$ action, we need some computation because in this case the
size of the determinant changes.
For $w(\Lambda_0) \leftrightarrow (I,J)$ such that $m_0=1$, we have
$s_0 w(\Lambda_0) \leftrightarrow (I \cup \{1\},J \cup \{0\})$.
Then 
$$
\tp_{s_0 w}=\det
\left[
\matrix{
\left(\tp_{p,q}\right)_{p \in I, q \in J}&\left(\tp_{1,q}\right)_{q \in J} \cr
\left(\tp_{p,0}\right)_{p \in I}&\tp_{1,0}
}
\right].
$$
We shall prove that this is equal to $s_0(\tp_w) {f_0 \over \alpha_0}$.

As is shown in eq.(1), we have
$$
s_0(\tp_{2,1}) \tp_{1,0}=\tp_{2.1} \tp_{1,0}-\tp_{2,0} \tp_{1,1}.
$$
Applying the actions $s_k$, $k \geq 2$ or $k \leq -2$ repeatedly, one get
$$
s_0(\tp_{p,q}) \tp_{1,0}=\tp_{p.q} \tp_{1,0}-\tp_{p,0} \tp_{1,q},
$$
for any $p \geq 2$ and $q \geq 1$.
By using this and $\tp_{1,0}=f_0/\alpha_0$, we get
$$
{f_0 \over \alpha_0} s_0\left(\tp_w(\Lambda_0)\right)=
\tp_{1,0}^{1-\vert I \vert}
\det \left[\tp_{p,q} \tp_{1,0} -\tp_{p,0} \tp_{1,q} \right]_{p \in I,q \in J},
$$
and this is indeed equal to $\tp_{s_0 w}(\Lambda_0)$, 
because of the identity in Lemma 5.
$\Box$

\vskip5mm
\noindent
{\it{Proof of the Theorem 2.}}
The proof for the general
$\tp_w(\Lambda_0)$ cases can be reduced to the single column cases by the
action of
$s_k$, ($k>0$) which preserves the size of the determinant.

For the single column case $\lambda=(1^{q+1})$, the desired formula is
\begin{equation}
\tp_{1,q}(\Lambda_0)=
\det
\left[
\matrix{
h_{1}^{(-q)}&h_{2}^{(-q)}&\cdots&h_{q+1}^{(-q)} \cr
1&h_{1}^{(-q+1)}&\cdots&h_{q}^{(-q+1)} \cr
&\ddots&\ddots&\vdots \cr
&1&h_{1}^{(-1)}&h_{2}^{(-1)} \cr
&&1&h_{1}^{(0)} \cr
}
\right].
\end{equation}
The case of $q=1$ follows from the relation in eq.(2).
More generaly, we have
\begin{equation}
\tp_{i,1}=\det
\left[
\matrix{
h_1^{(-1)}&h_{i+1}^{(-1)} \cr
1&h_{i}^{(0)} \cr
}
\right].
\end{equation}
To prove the case $q>1$, we will check the condition in Proposition 2,
\begin{equation}
{f_{-(q+1)} \over \alpha_{-(q+1)}} s_{-(q+1)}\big(\tp_{1,q}(\Lambda_0)\big)=
\tp_{1,q+1}(\Lambda_0).
\end{equation}
Using the eq.(5), the left hand side of this relation (6) can be written as 
the same determinant as $\tp_{1,q}(\Lambda_0)$ in (4),
with the first row $h_{i}^{(-q)}$ replaced by ${h'}_{i}^{(-q)}$ such as
$$
{h'}_{i}^{(-q)}=
{f_{-(q+1)} \over a_{-(q+1)}} s_{-(q+1)} \pi^{-q}(\tp_{i,0})=
\pi^{-q}(\tp_{i,1})=\det
\left[
\matrix{
h_{1}^{(-q-1)}&h_{i+1}^{(-q-1)} \cr
1&h_{i}^{(-q)} \cr
}
\right].
$$
Then the left hand side of the eq.(6) is equal to the right hand side
$$
\tp_{1,q+1}(\Lambda_0)=
\det
\left[
\matrix{
h_{1}^{(-q-1)}&h_{2}^{(-q-1)}&\cdots&h_{q+2}^{(-q-1)} \cr
1&h_{1}^{(-q)}&\cdots&h_{q+1}^{(-q)} \cr
&\ddots&\ddots&\vdots \cr
&1&h_{1}^{(-1)}&h_{2}^{(-1)} \cr
&&1&h_{1}^{(0)} \cr
}
\right],
$$
because of the identity
$$
\det
\left[
\matrix{
a_{00}\cr
a_{10}\cr
0\cr
\vdots\cr
0\cr
}
\left(
\matrix{
a_{0j}\cr
a_{1j}\cr
a_{2j}\cr
\vdots\cr
a_{nj}\cr
}
\right)_{1 \leq j \leq n}
\right]=
\det
\left[
\left(
\matrix{
\det
\left[
\matrix{
a_{00}&a_{0j}\cr
a_{10}&a_{1j}\cr
}
\right]
\cr
a_{2j}\cr
\vdots\cr
a_{nj}\cr
}
\right)_{1 \leq j \leq n}
\right].
$$
$\Box$

\section{Remarks}
\noindent
{\bf{Remark 1.}}
The reduction to the finite rank cases $A^{(1)}_{N-1}$ is given by
the $N$-reduced condition "$\pi^N=1$".
On the variables 
$\alpha,f, \tau$, this reduction is simply realized by the specialization
$$
\alpha_{i+N}=\alpha_i, \quad f_{i+N}=f_i, \quad \tau_{i+N}=\tau_i.
$$ 
By putting
$\bar{s_i}=\prod_{n \in \bz} s_{i+n N}$, $(i \in {\bz/N \bz})$,
the representation $W(A_{\infty})$ reduces to that of $W=W(A^{(1)}_{N-1})$ 
on the field ${\bf{C}}(\alpha_i;f_i;\tau_i, (i \in {\bz/N \bz}))$.
The $\tau$-functions for $A^{(1)}_{N-1}$ case is nothing but
the specialization of $\tau$ for $A_{\infty}$ and have the same 
determinant formulas.
Note that only the $N$-reduced Young diagrams can be generated by the
actions of $\bar{s_i}$.

\vskip5mm
\noindent
{\bf{Remark 2.}}
The determinant formulas of the $\tau$-function provide explicit
solutions for initial value problem
of the associated discrete dynamical system introduced in \cite{affine}.
The polynomiality of the $\tau$-functions and the multiplicative
formula of $f$-variables in terms of the $\tau$-functions
$$
w(f_i)={\phi_w(\Lambda_{i}) \phi_{ws_i}(\Lambda_i) \over
\phi_{w}(\Lambda_{i-1}) \phi_{w}(\Lambda_{i+1})},
$$
give a strong support for the singularity confinement property
which the discrete dynamical system expected to have.

\vskip5mm
\noindent
{\bf{Remark 3.}}
The representation of affine Weyl groups $W(A^{(1)}_l)$ has an interpretation
as the B\"acklund transformations for the Painlev\'e equations
$P_{IV}$ (for $l=2$), $P_{V}$ (for $l=3$) and their generalizations for
$l \geq 4$ \cite{higher}.
Under this interpretation, the polynomials $\phi_{w}(\Lambda_j)$ are the
far-reaching generalization of the "special polynomials" arising in Painlev\'e
equations in the sense of Umemura et. al. \cite{Um1}\cite{Um2}.
When specialized to certain "initial solutions", we obtain  explicit
determinant formulas for the Okamoto polynomials (for $P_{IV}$ \cite{PIV}), 
the Umemura polynomials (for $P_{V}$ \cite{PV}) and their generalizations. 

\vskip5mm
\noindent
{\bf{Remark 4.}}
The normalized functions $\tp_w(\Lambda_0)$ can be represented as the
minor determinants of the following frame $X$ of the universal Grassmann
manifold,
$$
X=
\left[
\matrix{
\ddots&\ddots&\ddots&\ddots&\ddots&\cdots
\cr &h_{0}^{(-2)}&h_{1}^{(-2)}&h_{2}^{(-2)}&h_{3}^{(-2)}&\cdots \cr
&&h_{0}^{(-1)}&h_{1}^{(-1)}&h_{2}^{(-1)}&\cdots \cr 
&&&h_{0}^{(0)}&h_{1}^{(0)}&\cdots \cr
}
\right].
$$
Where $h_{i}^{(j)}=\pi^j \big(\tp_{i,0}\big)$.
In this picture, the Weyl group $W$ is nothing but the Weyl group $W(GL_{\infty})$.
It would be interesting if the space of the initial values of Painlev\'e
equations can be realized as a natural sub-manifolds of 
the universal Grassmann manifold.

\vskip10mm
\noindent
{\bf {Acknowledgment.}}
I wish to thank Masatoshi Noumi for stimulating discussions.

\end{document}